\documentclass[12pt]{amsart}
\allowdisplaybreaks
\usepackage{amssymb}
\usepackage{graphicx}
\usepackage{fullpage}
\usepackage[pdfencoding=auto,colorlinks,linkcolor=,citecolor=]{hyperref}
\usepackage{cite}
\usepackage{dsfont}
\usepackage{eucal}
\usepackage[all,cmtip]{xy}

\newcommand{\gr}{\mathsf{gr}}

\newcommand{\Jen}{\mathsf{Jen}}

\newcommand{\bF}{\mathbb F}

\newcommand{\cS}{\mathcal S}

\newcommand{\cU}{\mathcal{U}}

\renewcommand{\le}{\leqslant}
\renewcommand{\ge}{\geqslant}

\makeindex
\numberwithin{equation}{section}

\theoremstyle{plain}

\newtheorem{theorem}[equation]{Theorem}

\theoremstyle{definition}

\theoremstyle{remark}

\author{David J. Benson} 
\address{Institute of Mathematics \\ 
University of Aberdeen \\ 
Aberdeen AB24 3UE \\ 
United Kingdom}
\email{d.j.benson@abdn.ac.uk}

\title{The socle of the group algebra of a finite $p$-group}



\begin{document}

\begin{abstract}
  Let $G$ be a finite $p$-group, and $\alpha$ an automorphism of the group
  algebra $\bF_pG$. Then $\alpha$ fixes the socle of $\bF_pG$
  pointwise. More generally, if $k$ is a field of characteristic $p$,
  and $\alpha$ is a $k$-algebra automorphism of $kG$, then $\alpha$
  induces a linear action on the dimension
  subquotients of the group, and the action on the socle is scalar
  multiplication by the $(p-1)$st power of the
  product of the determinants of this action.
  The scalar is thus an element of $(k^\times)^{p-1}$.
\end{abstract}

\maketitle

\section{Introduction}

Throughout this paper, we
let $G$ be a finite $p$-group and $k$ a field of characteristic $p$.
It is well known that the socle of $kG$ is one dimensional, spanned
by the sum of the group elements. Therefore given any algebra automorphism
$\alpha$ of $kG$, there is an element $\lambda\in k^\times$ such that
the action of $\alpha$ on the socle is given by multiplication by
$\lambda$.
We use the Jennings--Quillen theory, describing the associated graded
of the radical filtration to prove the following theorem.

\begin{theorem}\label{th:Rueping}
For any algebra automorphism $\alpha$ of $kG$, we have
$\lambda\in (k^\times)^{p-1}$. In particular if $k=\bF_p$ then
$\lambda=1$. 
\end{theorem}

The exact value of $\lambda$ can be computed as the $(p-1)$st
power of the determinant of a matrix given by the action on the
dimension subquotients of $G$.\bigskip

\noindent
{\bf Acknowledgement.} This paper grew out of a question posed by
Henrik R\"uping on the online forum {\sf MathOverflow}.

\section{The Jennings--Quillen theory}

In this section, we prove Theorem~\ref{th:Rueping}.
Recall that $G$ is a finite $p$-group and $k$ is a field of
characteristic $p$.
Jennings~\cite{Jennings:1941a} described the radical layers of $kG$
in terms of a certain descending central series of $G$, and
Quillen~\cite{Quillen:1968a} reinterpreted Jennings' theorem in terms
of restricted enveloping algebras. See also Section~3.14
of~\cite{Benson:1991a} and Section~V.7 of~\cite{Jacobson:1962a}
for background.

Let $J(kG)$ be the Jacobson radical of $kG$, spanned by the elements
$g-1$ for $g\in G$. The associated graded algebra is defined by
\[ \gr_*(kG)=\bigoplus_{n\ge 0}J^n(kG)/J^{n+1}(kG). \]
For $r\ge 1$, the dimension groups $F_r(G)$
are defined by
\[ F_r(G)=\{g\in G\mid g-1\in J^r(kG)\}. \]
Thus $F_1(G)=G$, $F_2(G)=\Phi(G)$, $[F_r(G),F_s(G)]\subseteq
F_{r+s}(G)$, and if $g\in F_r(G)$ then $g^p\in
F_{pr}(G)$. Furthermore, $F_r(G)$ is the most rapidly descending
central series with these properties. We define
\[ \Jen_*(G)=\bigoplus_{r\ge 1}k\otimes_{\bF_p}F_r(G)/F_{r+1}(G). \]
Then $\Jen_*(G)$ is a $p$-restricted Lie algebra over $k$ with Lie bracket
induced by taking commutators in $G$ and with the $p$-restriction map
$x \mapsto x^{[p]}$ induced by taking $p$th powers in $G$. As a
$p$-restricted Lie algebra, $\Jen_*(G)$ is generated by its degree one
elements. Let $\cU\Jen_*(G)$ be the restricted universal enveloping
algebra of $\Jen_*(G)$ over $k$. As an associative algebra,
$\cU\Jen_*(G)$ is generated by $\Jen_*(G)$, subject to relations given
by the commutators and $p$th powers. Note that multiplicative
commutators in $G$ become additive commutators in $\cU\Jen_*(G)$.

The Jennings--Quillen theorem states that there is a $k$-algebra
isomorphism
\begin{equation}\label{eq:JQ}
  \cU\Jen_*(G) \to \gr_*(kG)
\end{equation}
which, for any $r$ and any $g\in F_r(G)$, sends the image of $g$ in
$F_r(G)/F_{r+1}(G)$ to the image of $g-1$ in $\gr_*(kG)$.

Now let $|G|=p^m$. Then $\dim_k\Jen_*(G)=m$ and $\dim_k\cU\Jen_*(G)=p^m$.
Let
$x_1,\dots,x_m$ be a basis of homogeneous elements in
$\Jen_*(G)$, with $x_j$ the image of an element $g_j$ in the
corresponding dimension subgroup of $G$, for $1\le j\le m$.
We write $\bar x_j$ for the corresponding generators of $\cU\Jen_*(G)$.

Next, we introduce an increasing filtration on $\cU\Jen_*(G)$ as in
the proof of the Poincar\'e--Birkhoff--Witt theorem for Lie algebras. Let
$\cU^{(0)}\Jen_*(G)=k\cdot 1$ in degree zero,
$\cU^{(1)}\Jen_*(G)=k\cdot 1\oplus\Jen_*(G)$,
and for $i\ge 2$,
$\cU^{(i)}\Jen_*(G)=\Jen_*(G)\cU^{(i-1)}\Jen_*(G)$. Thus
$\cU^{(i)}\Jen_*(G)$ is the subspace of $\cU\Jen_*(G)$
spanned by words of length at most $i$ in the generators $\bar x_j$
of $\cU\Jen_*(G)$. Then
\begin{align*}
  \cU^{(i)}\Jen_*(G).\cU^{(j)}\Jen_*(G)&\subseteq\cU^{i+j}\Jen_*(G),\\
[\cU^{(i)}\Jen_*(G),\cU^{(j)}\Jen_*(G)]&\subseteq\cU^{i+j-1}\Jen_*(G),
\end{align*}
and if $x$ is in $\cU^{(1)}\Jen_*(G)$ then $x^p$ is also in
$\cU^{(1)}\Jen_*(G)$. Then the associated graded
\begin{equation}\label{eq:SJen}
  \cS\Jen_*(G)=\bigoplus_{i\ge
    0}\cU^{(i)}\Jen_*(G)/\cU^{(i-1)}\Jen_*(G)
\end{equation}
carries a double grading. The first degree keeps track of the
dimension subgroups, and the second degree comes from the filtration.
It is generated by $\cU^{(1)}\Jen_*(G)$, is commutative, and $p$th powers
of positive degree elements are zero. Writing $\bar{\bar x}_j$ for the
image of $\bar x_j$ in $\cS\Jen_*(G)$, the analogue of the
Poincar\'e--Birkhoff--Witt theorem for $p$-restricted Lie algebras
shows that
\[ \cS\Jen_*(G)\cong k[\bar{\bar x}_1,\dots,\bar{\bar x}_m]/
  (\bar{\bar x}_1^p,\dots,\bar{\bar x}_m^p), \]
is the truncated symmetric algebra on $\Jen_*(G)$.
Thus the elements $\bar x_1^{i_1}\dots\bar x_m^{i_m}$ form a basis
for $\cU\Jen_*(G)$, where $0\le i_j < p$ for $1\le j\le m$, and the
elements $\bar{\bar x}_1^{i_1}\dots\bar{\bar x}_m^{i_m}$ form a basis
for $\cS\Jen_*(G)$.

\section{Proof of the Theorem}

The socle
of $kG$ is one dimensional, spanned by the sum of the group elements,
\[ n=\sum_{g\in G}g=\prod_{j=1}^m(g_j^{p-1}+g_j^{p-2}+\dots+g_j+1)
  =\prod_{j=1}^m(g_j-1)^{p-1} \in kG. \]
The image of this in the restricted universal enveloping algebra $\cU\Jen_*(G)$, using the
isomorphisms~\eqref{eq:JQ}, is
$\bar n = \prod_{j=1}^m \bar x_j^{p-1}$. The image in
$\cS\Jen_*(G)$, using~\eqref{eq:SJen}, is
$\bar{\bar n}=\prod_{j=1}^m\bar{\bar x}_j^{p-1}$.

Now let $\alpha$ be a $k$-algebra automorphism of the group algebra
$kG$. Then $\alpha$ induces an automorphism $\bar\alpha$ of
$\cU\Jen_*(G)$ and an automorphism $\bar{\bar\alpha}$ of
$\cS\Jen^*(G)$. Since the socle of $kG$ is one dimensional,
the action of $\alpha$ on $n$ is multiplication by some scalar
$\lambda$, so $\alpha(n)=\lambda n$ with $\lambda\in k^\times$.
Taking associated graded with respect to the radical filtration and using
the isomorphism~\ref{eq:JQ}, we obtain $\bar\alpha(\bar n)=\lambda\bar
n$. Then, taking the associated graded with respect to the
Poincar\'e--Birkhoff--Witt filtration using~\eqref{eq:SJen}, we have
\begin{equation}\label{eq:barbar}
  \bar{\bar\alpha}(\bar{\bar n})=\lambda\bar{\bar n}.
\end{equation}

To compute the value of $\lambda$, we look at the action of $\alpha$
on $\cS\Jen_*(G)$. This preserves the double grading, and so it is
given by an $m\times m$ invertible matrix $A$ on the generators
$\bar{\bar x}_1,\dots,\bar{\bar x}_m$. Denote by $\det(A)$ the
determinant of the matrix $A$. Then
\begin{equation}\label{eq:lambda}
  \bar{\bar\alpha}(\bar{\bar n}) =
  \prod_{j=1}^m \bar{\bar\alpha}(\bar{\bar x}_j)^{p-1}
  =\det(A)^{p-1}\prod_{j=1}^m\bar{\bar x}_j^{p-1}
  =\det(A)^{p-1}\bar{\bar n}.
\end{equation}
To check the second equality, it suffices to check it on generators
for the general linear group. Taking the generators to be the
elementary matrices and the diagonal matrices, the check is easy.

It follows from~\eqref{eq:barbar} and~\eqref{eq:lambda} that
\[ \lambda = \det(A)^{p-1}\in (k^\times)^{p-1}, \]
and this completes the proof of Theorem~\ref{th:Rueping}.

\bibliographystyle{amsplain}
\bibliography{../repcoh}

\end{document}